\documentclass[12pt,thmsa,emstex]{article}
\usepackage[left=3cm,right=3cm,top=3cm,bottom=3cm]{geometry}
\usepackage{amsfonts}
\usepackage{t1enc}
\usepackage{amssymb}
\usepackage{epsfig}
\usepackage[french,english]{babel}
\usepackage{amsmath}
\usepackage{graphics}
\usepackage{layout}
\usepackage{latexsym}
\usepackage{color}
\usepackage{authblk}

\def\wdt {\widetilde}

\newtheorem{theorem}{Theorem}[section]
\newtheorem{lemma}[theorem]{Lemma}
\newtheorem{corollary}[theorem]{Corollary}
\newtheorem{proposition}[theorem]{Proposition}
\newtheorem{example}[theorem]{Example}
\newtheorem{remark}[theorem]{Remark}

\newtheorem{hypothesis}[theorem]{Assumption}

\newtheorem{definition}[theorem]{Definition}
\newtheorem{exercice}[theorem]{Exercice}

\def\bit{\begin{itemize}}
	
	\def\eit{\end{itemize}}

\reversemarginpar   

\def\bc{\begin{center}}
	
	\def\ec{\end{center}}

\def\bthm{\begin{theorem}}
	
	\def\ethm{\end{theorem}}

\def\bcor{\begin{corollary}}
	
	\def\ecor{\end{corollary}}

\def\bprop{\begin{proposition}}
	
	\def\eprop{\end{proposition}}

\def\blem{\begin{lemma}}
	
	\def\elem{\end{lemma}}

\def\bex{\begin{example}}
	
	\def\eex{\end{example}}

\def\bexo{\begin{exercice} \rm }
	
	\def\eexo{\end{exercice} }

\def\brem{\begin{remark}}
	
	\def\erem{\end{remark}}

\def\prf{{\bf Proof }}

\def\bdes{\begin{description}}
	
	\def\edes{\end{description}}

\def\ita{\item[(a)]}

\def\itb{\item[(b)]}

\def\itc{\item[(c)]}

\def\iti{\item[(i)]}

\def\itii{\item[(ii)]}

\def\beq{\begin{equation}}

\def\eeq{\end{equation}}

\def\ben{\begin{enumerate}}
	
	\def\een{\end{enumerate}}

\def\beqar{\begin{eqnarray}}

\def\eeqar{\end{eqnarray}}

\def\beqarr{\begin{eqnarray*}}
	
	\def\eeqarr{\end{eqnarray*}}

\def\QED{\hfill\ensuremath{\square}\\[2ex]}

\def\prf{{\bf Proof }\hspace{.1in}}

\def\L{{\cal L}}
\def\D{{\cal D}}

\def\E{{\mathsf E}}

\def\M{{\cal M}}

\def\RR{{\mathbb R}}  

\def\EE{{\mathbb E}}
\def\E{{\mathcal E}}
\def\PP{{\mathbb P}}

\newcommand{\bdelta}{\boldsymbol{\delta}}
\def\eps{\varepsilon}

\newcommand{\Pcal}{\mathcal{P}}

\newcommand{\Fcal}{\mathcal{F}}

\def\1{{\rm 1\mskip-4.4mu l}}
\begin{document}
	\title{A method to deal with the critical case in stochastic population dynamics}
	\author[1]{Dang H. Nguyen }
	
	\author[2]{Edouard Strickler}
	
	\affil[1]{ University of Alabama, Tuscaloosa, USA}
	\affil[2]{ Université de Lorraine, CNRS, Inria, IECL, Nancy, France }
	\maketitle
	\begin{abstract}
		In numerous papers, the behaviour of  stochastic population models is investigated through the sign of a real quantity which is the growth rate of the population near the extinction set. In many cases, it is proven that when this growth rate is positive, the process is persistent in the long run, while if it is negative, the process converges to extinction. However, the critical case when the growth rate is null is rarely treated. The aim of this paper is to provide a method that can be applied in many situations to prove that in the critical case, the process congerves in temporal average to extinction. A number of applications are given, for Stochastic Differential Equations and Piecewise Deterministic Markov Processes modelling prey-predator, epidemilogical or structured population dynamics. 
		
	\end{abstract}
	\section{Introduction}
	
	Much effort in population biology has been devoted to understanding under what conditions interacting populations, whether they be viruses, plants, or animals, coexist or go extinct.
	The variation of environmental factors such as temperature, precipitation and humidity inherently affect the growth rates of the species. There is rich literature showing that 
	the interplay of biotic interactions and environmental fluctuations  can facilitate or suppress the persistence of species or disease prevalence; see \cite{gillespie1978effects,chesson1981environmental, abrams1998apparent, chesson2008interaction, BL16} and the references therein. There has been intensive attention paid to modeling and analysis of ecological and epidemiological models under environmental stochasticity.
	
	In \cite{SBA11}, a condition for coexistence was
	given, which requires  a certain weighted combination of populations’ invasion rates to be positive for any invariant measures associated with a subcollection of populations. 
	The results were then improved and generalized to a very general setting by Michel Benaïm in \cite{B18},
	where the concept of \emph{$H$-persistence} was coined and developed.
	With the same idea, \cite{HN16} provided conditions for both extinction and persistence in a setting of stochastic differential equations.
	The long-term properties of some specific models are also classified in \cite{DNDY16, DNY16, BL16, NY17, HS17, BS17, GPS19}.
	For many models, the conditions in the aforesaid references for extinction and persistence
	of a species in an interacting populations are determined by a threshold $\Lambda$ whose sign indicates whether the species will be persistent or extinct.
	Namely, the  result obtained is that
	if $\Lambda>0$ the species persists
	and if $\Lambda<0$, extinction will happen.
	While this kind of results are sharp in the sense that they leave only critical cases $(\Lambda=0)$ untreated,
	it is of great interest to discover the dynamics of the systems in critical cases.
	Similar to the case of an equilibrium of a deterministic dynamical system whose maximum eigenvalue is 0, treating the critical cases of stochastic systems is, in general, extremely difficult. However, populations models often exhibit
	some certain monotone properties that 
	can be utilized to handle critical cases.
	This paper provides some methods for treating the critical cases
	of population dynamics under certain conditions. It is partially inspired by the work of the first author \cite{NY17} where the critical case is treated for a stochastic chemostat dynamic modeled by a switching diffusion. 
	
	The rest of the paper is organized as follows.
	In Section 2, we formulate the model in the general setting of \cite{B18} and
	give a general condition for extinction in average of stochastic populations in a critical case.
	Section 3 is devoted to the analysis of a number of specific models in critical cases.
	Different techniques are introduced so that
	the general result in Section 2 become applicable for those models.
	
	\section{Notations and results}
	
	Before to give our result,  we present the very general framework of \cite{B18} for stochastic persistence and extinction. Let $(X_t)_{ t \geq 0}$ be a cadlag Markov Process on a locally compact Polish metric space $(\M,d)$. For a distribution $\nu$ on $\M$, we set, as usual, $\PP_{\nu}$ for the law of the process $X$ with initial distribution $\nu$ and $\EE_{\nu}$ for the associated expectation. If $\nu = \delta_x$ for some $x \in \M$, we write $\PP_x$ for $\PP_{\delta_x}$. We denote by $(P_t)_{t \geq 0}$ the semigroup of $X$ acting on bounded measurable function $f : \M \to \RR$ as 
	$$ 
	P_t f(x) = \EE_x \left( f(Z_t) \right).
	$$
	An invariant distribution for the process $Z$ is a probability $\mu$ such that $\mu P_t = \mu$ for all $t \geq 0$. We let $\Pcal_{inv}$ denote the set of all the invariant distributions of $X$ and for $N \subset \M$, let  $\Pcal_{inv}(N)$ and $\Pcal_{erg}(N)$ denote the (possibly empty) sets of invariant probability measures and ergodic invariant probability measures, respectiveley,  giving mass 1 to the set $N$. The following is the standing assumption:
	
	\begin{hypothesis}
		\label{hyp:extinct}
		There exists a non empty closed set ${\M}_0 \subset \M$  called the {\em extinction set} which is invariant under $(P_t)_{t \geq 0}.$ That is, for all $t \geq 0$,
		$$P_t \1_{\M_0} = \1_{\M_0}.$$
	\end{hypothesis} We set $$\M_+ = \M \setminus \M_0.$$
	The two following assumptions are taken from \cite{B18}.
	
	\begin{hypothesis}
		\label{hyp:feller}
		The semigroup $(P_t)_{t \geq 0}$ is \emph{$C_b$ - Feller}, meaning that for all continuous bounded function $f : \M \to \RR$,  $ (t,x) \mapsto P_t f(x)$ is a  continuous function. 
	\end{hypothesis}
	We let $\L$ denotes the infinitesimal generator of $P_t$ on the space $C_b(\M)$ of continuous bounded functions, defined for $f \in \mathcal{D}(\L)$ by
	\[
	\L f(x) = \lim_{t \to 0} \frac{P_t f(x) - f(x)}{t},
	\]
	where the domain is the set of functions such that the above convergence holds pointwise, with the additionnal property that $ \sup_{0 < t \leq 1} \| t^{-1}( P_t f - f) - \L f\| < + \infty$. 
	We also let $\D^2(\L)$ be the set of functions such that both $f$ and $f^2$ lie in $\D (\L)$, and we define the carré du champ operator on $\D^2 (\L)$ by 
	\[
	\Gamma f = \L f^2 - 2 f \L f. 
	\]
	For all $t > 0$, we let $\Pi_t$ denote the empirical occupation measure of the process $X$ up to time $t$. This is the random probability measure defined on $\M$ by
	$$
	\Pi_t = \frac{1}{t} \int_0^t \delta_{X_s} ds.
	$$
	When we want to emphasis the starting point, we set $\Pi_t^x$ for the empirical occupation measure whenever $X_0=x$ almost surely.
	\begin{hypothesis}
		For all $x \in M$, the sequence $\{ \Pi_t^x, t \geq 0\}$ is almost surely tight.
	\end{hypothesis}
	As it is proved in \cite[Theorem 2.1 ]{B18}, a sufficient condition for the tightness of the sequence of the empirical occupation measures is the existence of a suitable Lyapunov function, as defined in the following assumption. Recall that a map $f: \M \to \RR_+$ is said to be \textit{proper} if for all $R > 0$, the sublevel set $\{ f \leq R\}$ is compact in $\M$. 
	
	\begin{hypothesis}
		\label{hyp:lyapinfty}
		There exist  continuous proper  maps $W, \tilde W: \M \mapsto \RR_+$ and a continuous map $LW : \M \mapsto \RR$ enjoying the following properties :
		\bdes
		\ita For all compact $K \subset \M$ there exists $W_K \in \D^2$ with  $W|_K = W_K|_K$ and  $(\L W_K)|_K = LW|_K;$
		\itb For all $x \in M$, $\sup_{\{t\geq 0, K : K \subset \M, K \mbox{ compact } \}} P_t \Gamma(V_K)(x)   < \infty;$
		\itc $LW \leq - \tilde{W} + C$.
		\edes
	\end{hypothesis}
	The latter assumption also implies that all weak-limit point of the sequence $(\Pi_t)_{t > 0}$ are almost surely in $\Pcal_{inv}(\M)$ (see \cite[Theorem 2.1 ]{B18}).
	
	The next assumption ensures the existence of a Lyapunov function near the boundary $\M_0$ :
	\begin{hypothesis}
		\label{hyp:VH}
		There exist  continuous maps $V : \M_+ \mapsto \RR_+$ and $H :\M \mapsto \RR$ enjoying the following properties :
		\bdes
		\ita For all compact $K \subset \M_+$ there exists $V_K \in \D^2$ with  $V|_K = V_K|_K$ and  $(\L V_K)|_K = H|_K;$
		\itb For all $x \in \M$, $\sup_{\{K : K \subset \M, K \mbox{ compact }; \: t \geq 0 \}} P_t \Gamma(V_K)(x)   < \infty;$
		\itc The map $\frac{\tilde{W}}{1+|H|}$ is proper.
		\edes
	\end{hypothesis}
	
	From this assumption, it is possible to define the $H$ - exponent of $X$ as in \cite[Definition 4.2]{B18}.
	
	\begin{definition}
		For $V$ and $H$ as in Assumption , we set 
		
		$$
		\Lambda^-(H) = - \sup \{ \mu H, \: \mu \in \Pcal_{erg}(\M_0) \}
		$$
		and
		$$
		\Lambda^+(H) = - \inf \{ \mu H, \: \mu \in \Pcal_{erg}(\M_0) \}.
		$$
		We say that $X$ is $H$ \emph{- persistent} if $\Lambda^-(H) > 0$ and that $X$ is \emph{$H$ - nonpersistent} if $\Lambda^+(H) < 0$.
	\end{definition}
	
	The main results in \cite{B18} could be summed up as follows.  If $\Lambda^-(H) > 0$, then $\Pcal_{inv}(M_+)$ is non empty and the family $\{\Pi_t, t \geq 0\}$ is tight in $\M_+$. Furthermore, the process $X$ is \textit{stochastically persistent } (see \cite{Sch12})) in the sense that, for all $\eps > 0$, there exists a compact subset $K$ of $\M_+$ such that, for all $x \in \M_+$, 
	\[
	\PP_x\left( \liminf \Pi_t(K) \geq 1 - \eps \right)=1.
	\]
	On the contratry, when $\Lambda^+(H) < 0$,  $X_t$ converges to $\M_0$ exponentially fast (this is not yet proven in \cite{B18}, but one can look at the thesis of the secound author \cite[Section 1.3]{these} for a proof in the special case where $\M_0$ is compact, relying on the proof made in \cite{BL16}). However, the critical case where $\Lambda^+(H) = 0$ is not investigated. It is  known from the deterministic case that in general, the information  that $\Lambda^+(H) = 0$ is not sufficient to conclude on the long term behaviour of the process (one can think to the stability of an equilibrium point for a dynamical system, when the Jacobian matrix of the vector field at that point has eigenvalues with null real part).
	
	We now state the result of this note, which follows from a basic argument :
	
	\bprop
	\label{prop}
	Assume that if $\Pcal_{inv}(\M_+)$ is non empty, then  there exists $\mu \in \Pcal_{inv}(\M_+)$ and $\pi \in \Pcal_{inv}(\M_0)$ such that
	
	\begin{equation}
	\label{eq:mugreaterpi1}
	\mu H > \pi H.
	\end{equation}
	Then $\Lambda^+(H) > 0$.
	
	\eprop
	
	\prf
	Assume  that $\Pcal_{inv}(\M_+)$ is nonempty. Let $\mu \in \Pcal_{inv}(\M_+)$ satisfying \eqref{eq:mugreaterpi1} for some $\pi \in \Pcal_{inv}(\M_0)$, then  $\mu H > - \Lambda^+(H)$. By 
	\cite[Lemma 7.5]{B18}, since $\mu \in \Pcal_{inv}(\M_+)$, we must have $\mu H = 0$ (note that the proof of this fact in \cite{B18} does not require the process to be $H$ - persistent.) This proves that $\Lambda^+(H) > 0$.
	\QED
	We get the following immediate corollary 
	
	\bcor
	\label{cor:lambda0}
	Assume that the hypothesis in Proposition \ref{prop} holds.
	If $\Lambda^+(H) = 0$,  $\Pcal_{inv}(\M_+)$ is empty and all weak-* limit point of $\Pi_t$ lie almost surely in $\Pcal_{inv}(\M_0)$. In particular, if  $\Pcal_{inv}(\M_0) = \{ \pi \}$, then for all bounded continuous function $f : \M \to \RR$, 
	\begin{equation}
	\label{eq:cvergoM0}
	\lim_{t \to + \infty} \frac{1}{t} \int_0^t f(X_s) ds = \pi f.
	\end{equation}
	\ecor
	
	\brem
	\label{rem:proper}
	Actually, one can prove that \eqref{eq:cvergoM0} holds for all $f : \M \to \RR$ such that the map $\frac{W}{1+|f|}$ is proper, where $W$ satisfy Assumption \ref{hyp:lyapinfty} (see \cite[Lemma 9.1]{B18}).
	\erem
	Thus, the idea is that if $H$ is strictly bigger on $\M_+$ than on $\M_0$ and if $\Lambda^+(H)= 0$, then the process goes in average to extinction. Rather than giving abstract conditions ensuring that \eqref{eq:mugreaterpi1} holds, we provide in the next sections five examples on which we prove \eqref{eq:mugreaterpi1} with different methods, that can be easily reproduce for other models. 
	\section{Applications}
	In this section, we prove that the results of the previous sections apply to five models. The four first examples come from the literature, where the case $\Lambda = 0$ has not be treated. The last example is new.
	
	\subsection{SIR model with switching}
	\label{sub:SIR}
	In this section, we apply our method to a SIRS model with random switching that was studied in \cite{li17}. We first describe the process. Let $N$ be a positive integer, and set $\E=\{1,\ldots,N\}$. For $k \in \E=\{1,\ldots,N\}$ let $F^k$ be the vector field defined on $\RR^3$ by:
	
	\begin{equation}\label{e:SIRS}
	F^k(S,I,R) =  \begin{pmatrix}
	\Lambda- \mu S + \lambda_k R - \beta_k S G_k(I)\\
	\beta_k S G_k(I)-(\mu+\alpha_k+\delta_k)I\\
	\delta_k I -(\mu + \lambda_k)R
	\end{pmatrix},
	\end{equation} 
	where $G_k$ is a regular function such that $G_k(0)=0$.  The reader is referred to \cite{li17} for the epidemiological interpretation of the different constants. Let $(\alpha_t)_{t \geq 0}$ be a irreducible Markov chain on $\E$. We denote by $p=(p_1, \ldots, p_N)$ its unique invariant probability measure. We consider the process $(Z_t)_{t \geq 0} = (X_t, \alpha_t)_{t \geq 0}$, with $X_t = (S_t, I_t, R_t) \in \RR_+^3$ evolving according to
	\begin{equation}
	\frac{d X_t}{dt} = F^{\alpha_t}(X_t).
	\end{equation}   
	The process $Z$ is a\textit{ Piecewise Deterministic Markov Process} (PDMP) as introduced in \cite{Dav84}, and belongs to the more specific class of PDMPs recently studied in \cite{BH12} and \cite{BMZIHP} (see also \cite{BL16}, \cite{HS17}, \cite{BS17} and \cite{GPS19} for PDMP model in ecology or epidemiology). 
	\brem
	In \cite{li17}, $\beta$ is the only parameter allowed to depend on $k$. The general case where the other constants and the function $G$ can depend on $k$ has been treated in \cite{S18}. 
	\erem
	We make the following assumptions, that are taken from \cite{li17} :
	
	\begin{hypothesis}\
		\label{hyp:SIR}
		\bdes
		\iti For all $k$, $G_k : \RR_+ \to \RR_+$ is $C^2$, with $G_k(0)=0$ and $0 < G_k(I) \leq G_k'(0)I$ for $I > 0$;
		\itii For all $k$, if $\beta_k \frac{\Lambda}{\mu} G_k'(0) - (\mu + \alpha_k + \delta_k) > 0$, then $F^k$ admits an equilibrium point $x^* \in \M_+$ which is accessible from $\M_+$.
		\edes
	\end{hypothesis}
	We consider the process on the space $\M:=K \times \E$, where $K = \{ (s,i,r) \in \RR_+^3 \: s + i + r \leq \frac{
		\Lambda}{\mu}\}.$ The set $K_0 = \{ I = 0 \}$ is invariant for the $F^k$ thus the set $\M_0 = K_0 \times E$ is invariant for $Z$. On this set, it is not hard to check that $X$ converge almost surely to $(S^*,0,0)$, where $S^* = \frac{\Lambda}{\mu}$. Thus, the unique invariant probability measure of $Z$ on $\M_0$ is $\delta^* \otimes p$, where $\delta^*$ is the Dirac mass at $(S^*,0,0)$. Consider the function
	$V : \M_+ \times \E \to \RR_+$ given by 
	$$V(s,i,r,k) = \log \frac\Lambda\mu-\log i\,\,\, \text{ for all } (s,i,r,k) \in \M_+ \times E.$$
	Define also the function $H : \M \times E \to \RR$ by $H(s,i,r,k) = (\mu + \alpha_k + \delta_k  - \beta_k s \tilde{G}_k(i))$ where $\tilde{G}_k$ is given by :
	$$
	\tilde{G}_k(i) = \begin{cases}
	\frac{G_k(i)}{i} \quad \mbox{if $i \neq 0$}\\
	G'(0) \quad \mbox{if $i = 0$}.
	\end{cases}
	$$
	It is not hard to check that $V$ and $H$ satisfy assumption \ref{hyp:VH}. Moreover, we have for $\pi = \delta^* \otimes p$, 
	
	$$
	\pi H = \sum_{k \in \E} p_k \left(\mu + \alpha_k + \delta_k  - \beta_k \frac{\Lambda}{\mu} G_k'(0)\right).
	$$
	for $k \in \E=\{1,\ldots,N\}$, As in \cite{li17}, we set  $$ R_0 = \frac{\sum_k p_k \beta_k \frac{\Lambda}{\mu} G'(0)}{\sum_k p_k (\mu + \alpha_k + \delta_k)}.$$ Note that $R_0 < 1$ (respectively $R_0>1$, $R_0 = 1$) if and only if $\pi H > 0$ (resp. $\pi H < 0$, $\pi H = 0$). The behaviour of the process when $R_0<1$ or $R_0 >1$ is studied in \cite{li17} (see also \cite{S18} for an alternative and more general proof). With our method, one can prove the following :
	\bprop
	Assume that $R_0 = 1$. Then, for all $(s,i,r,k) \in \M$, $\PP_{(s,i,r,k)}$ - almost surely,
	$$
	\lim_{t \to \infty} \frac{1}{t}\int_0^t S_u d u = S^*,
	$$
	and 
	$$
	\lim_{t \to \infty} \frac{1}{t}\int_0^t (I_u + R_u)  d u = 0.
	$$
	\eprop
	
	\prf We show that when $\Pcal_{inv}(\M_+)$ is nonempty, then for all $\mu^* \in \Pcal_{inv}(\M_+)$, one has 
	$$
	\mu^* H > \pi H.
	$$
	For convenience, we write $C_k$ for $\mu + \alpha_k + \delta_k$. By Assumption \ref{hyp:SIR}, we have 
	$$
	H(s,i,r,k) \geq C_k - \beta_k G'(0) s,
	$$
	and thus 
	$$
	\mu^* H \geq \sum_{k \in E} p_k C_k - \sum_{k \in E} \beta_k G_k'(0) \int_{M_+} s d \mu^*_k(s,i,r),
	$$
	where $\mu_k^*$ is the measure of total mass $p_k$ defined on $\M$ by $\mu^*_k(A) = \mu^*(A \times \{k\})$. Note that as $i > 0$ on $\M_+$ and that for $(s,i,r) \in \M$, $s+i+r \leq S^*$, then for all $(s,i,r) \in \M_+$, $s < S^*$. In particular, 
	
	$$
	\int_{\M_+} s d \mu^*_k(s,i,r) < p_k S^*,
	$$
	which yields 
	
	$$
	\mu^* H > \sum_{k \in \E} p_k C_k - \sum_{k \in \E} p_k \beta_k G_k'(0) S^* = \pi H.
	$$
	This proves by Corollary \ref{cor:lambda0} that if $R_0 = 1$, then $\Pcal_{inv}(\M_+)$ is empty and for all bounded measurable function $f : \M \times E \to \RR$,
	
	$$
	\lim_{t \to \infty} \frac{1}{t}\int_0^t f(S_u,I_u,R_u,r_u) d u = \sum_k p_k f(S^*,0,0,k).
	$$
	\QED
	
	\subsection{Stochastic Rosenzweig - MacArthur}
	\label{sub:RMA}
	This example is taken from \cite[Section 5.2]{B18}. We consider the following Stochastic Differential Equation (SDE) :
	
	\begin{equation}
	\begin{cases}
	d X_t = X_t \left( 1 - \frac{X_t}{K} - \frac{Y_t}{1+X_t} \right) dt + \eps X_td B_t\\
	d Y_t = Y_t \left( - \alpha + \frac{X_t}{1+X_t} \right)dt
	\end{cases}
	\end{equation}
	In this case, $\M = \RR_+^2 := \{(x,y) \in \RR^2 \: : x, y \geq 0 \}$. It is proven in \cite[Theorems 5.1 and 5.5]{B18} that Assumption \ref{hyp:lyapinfty} is satisfied with $W(x,y) = (x+y)^2$ and $\tilde{W} = (1 + C) W$, where $C$ is some constant. We set $\M_0^{\mathrm{x}} = \{ (x,y) \in \M \: : x = 0 \}$, $\M_0^{\mathrm{y}} = \{ (x,y) \in \M \: : y = 0 \}$ and $\M_0 = \M_0^{\mathrm{x}} \cup \M_0^{\mathrm{y}}$. We also let $\M_+^{\mathrm{x}} = \M \setminus \M_0^{\mathrm{x}}$,  $\M_+^{\mathrm{y}} = \M \setminus \M_0^{\mathrm{y}}$ and $\M_+ = \M \setminus \M_0$. We also define the invasion rate of species $x$ and $y$, respectively, as
	\[
	\lambda_1(x,y) = \left( 1 - \frac{x}{K} - \frac{y}{1+x} \right) - \frac{\eps^2}{2}
	\]
	and
	\[
	\lambda_2(x,y) =  - \alpha + \frac{x}{1+x}. 
	\]
	By \cite[Theorem 5.5]{B18}, if $\eps ^ 2 > 2$, then for any initial condition, one has $(X_t, Y_t) \to 0$ has $t \to \infty$. Thus, we assume now that $\eps^2 < 2$. In that case, the process is $H$ - persistent with respect to $\M_0^{\mathrm{x}}$. Indeed, in that situation, $\Pcal_{erg}(\M_0^{\mathrm{x}}) = \{ \delta_0 \}$, where $\delta_0$ is the Dirac mass at $0$ and $\delta_0 \lambda_1 = 1 - \frac{\eps^2}{2} > 0$. Hence, condition of \cite[Theorem 5.1 (ii)]{B18} is satisfied. In particular, every limit point of $(\Pi_t)_{t \geq 0}$ lies almost surely in $\Pcal_{inv}(\M_+^{\mathrm{x}})$. Moreover, on $\M_+^{\mathrm{x}} \cap \M_0^{\mathrm{y}}$, the process admits a unique invariant probability measure denoted by $\mu_{\mathrm{x}}$ ( see \cite[Section 5.2]{B18}).
	
	It is easily seen that $\Pcal_{erg}(\M_0) = \{ \delta_0, \mu_{\mathrm{x}} \}$ . We set 
	\[
	\Lambda(\eps, K, \alpha) = \mu_{\mathrm{x}}(\lambda_2) = \int_0^{+\infty} \frac{x}{1+x} d\mu_{\mathrm{x}}(x) - \alpha.
	\]
	By \cite[Theorem 5.5]{B18}, if $\Lambda(\eps, K, \alpha) > 0$, then the process is stochatistically persistent with respect to $\M_0$ and admits a unique invariant probability measure $\mu^*$ on $\M_+$, while if $\Lambda(\eps, K, \alpha) < 0$, $Y_t$ converges to $0$. We now prove the following proposition for the critical case :
	
	\bprop
	If $\Lambda(\eps, K, \alpha) = 0$, then for all $(x,y) \in M_+$, one has $\PP_{(x,y)}$ - almost surely, 
	\[
	\lim_{T \to \infty} \frac{1}{T} \int_0^T Y_s ds =0
	\]
	and
	\[
	\lim_{T \to \infty} \frac{1}{T} \int_0^T X_s ds = \int_0^{+ \infty} x d \mu_{\mathrm{x}} (x) =  K \left( 1 - \frac{\eps^2}{2} \right).
	\]
	\eprop

	\prf We prove that if $\Pcal_{inv}(\M_+^{\mathrm{y}})$ is non-empty, then for all $\mu^* \in \Pcal_{inv}(\M_+^{\mathrm{y}})$, one has $\mu^* H > \mu_{\mathrm{x}} H$, where 
	\[
	H(x,y) = H_1(x,y)-\lambda_2(x,y),
	\]
	with
	$$H_1(x,y)=\frac{1}{1+x+y}\left(x-\alpha y -\frac{x^2}K\right)-\frac{\eps^2x^2}{2(1+x+y)}.$$
	We set, for $(x,y) \in \M_+$, $V(x,y)=\log(1+x+y)-\log x$. We  can see that $(V,H)$
	satisfy Assumption \ref{hyp:VH}. Moreover, we have $\L [\log(1+x+y)]=H_1(x,y)$,
	then by \cite[Remark 19]{B18}, we must have
	$\nu H_1=0$ for any $\nu\in\Pcal_{inv}(\M)$. As a result,
	$$\nu H=-\nu\lambda_2\text{ for any }\,\nu\in\Pcal_{inv}(\M).$$
	\brem
	In the framework of \cite{B18}, it would have been natural to take for $V$ any function coinciding with $ -\log x$ for $x$ small enough, so that $H = - \lambda_2$ near $\M_0$, because it is sufficient to know $H$ on the boundary $\M_0$. However, to apply our method, it is required to compare $\pi H$ and $\mu H$ for $\mu \in \Pcal_{inv}(\M_+)$, thus it is necessary to know $H$ on the whole $\M_+$. Thus the idea is to take $V= V_1 + V_2$ and $H = H_1 + H_2$, with $V_2(x,y) = - \log x$, $H_2 = - \lambda_2$, $V_1$ defined on all $\M$  so that $V$ is nonnegative,  $\L V_1 = H_1$ and  $\nu H_1 = 0$ for all $\nu \in \Pcal_{inv}(\M)$ (see \cite[Remarks 11 and 19, and Proposition 4.13]{B18}). We use a similar trick in Subsection \ref{sub:patchy}.
	\erem
	
	To continue the proof, note that on $\M_0^{\mathrm{x}} \cap \M_+^{\mathrm{y}}$, $Y_t$ converges exponentially fast to $0$. Thus, it holds that $\Pcal_{inv}(\M_+^{\mathrm{y}}) = \Pcal_{inv}(\M_+)$. Moreover, by Theorem 5.5 in \cite{B18}, if  $\Pcal_{inv}(\M_+)$ is nonempty, it reduces to a unique element, that we denote by $\mu^*$, and $\mu^*$ has a positive density with respect to the Lebesgue measure. This implies by Birkhoff's ergodic theorem that for all $(x,y) \in \M_+$, 
	\[
	\mu^*\lambda_2 = \lim_{T \to \infty} \frac{1}{T} \int_0^T \lambda_2 (X_s, Y_s) ds, \quad \PP_{x,y} - \mbox{almost surely}.
	\]
	We let $\hat X$ be the solution of the reduced system on $\M_0^{\mathrm{y}}$. That is, 
	\begin{equation}
	d \hat X_t = \hat X_t \left( 1 - \frac{\hat X_t}{K} \right) dt + \eps d B_t.
	\end{equation}
	By the comparison theorem, if $X_0 = \hat X_0$, then $X_t \leq \hat X_t$ for all $t \geq 0$. The idea is now to write 
	\[
	\mu^*H  =-\mu^*\lambda_2= -\lim_{T \to \infty} \frac{1}{T} \int_0^T \lambda_2(\hat X_s, 0) ds 
	- \lim_{T \to \infty} \frac{1}{T} \int_0^T \left( \lambda_2(X_s, Y_s) - \lambda_2(\hat{X}_s, 0)\right) ds
	\]
	and to prove that the first term is $\mu_{\mathrm{x}} H$ and the second one is positive. 
	
	By \cite[Theorem 5.1 (i)]{B18}, we have $\mu_{\mathrm{x}}(\lambda_1) = 0$. Moroever, the process $\hat X$ on $\M_0^{\mathrm{y}}$ is persistent with respect to $\M_0^{\mathrm{x}} \cap \M_0^{\mathrm{y}}$. Thus, for all $x > 0$, one has 
	\[
	\lim_{T \to \infty} \frac{1}{T} \int_0^T \lambda_2(\hat X_s, 0) ds  = \mu_{\mathrm{x}} \lambda_2
	\]
	and
	\[
	\lim_{t \to \infty} \frac{1}{T} \int_0^T \lambda_1(\hat{X}_s,0) ds = \mu_{\mathrm{x}}(\lambda_1) = 0,
	\]  
	which gives
	\[
	\lim_{t \to \infty} \frac{1}{T} \int_0^T \hat X_s = K \left( 1 - \frac{\eps^2}{2} \right).
	\]
	On the other hand, since $(X,Y)$ is persistent with respect to $\M_0^{\mathrm{x}}$, one has 
	\[
	\lim_{t \to \infty}  \frac{1}{T} \int_0^T \lambda_1(X_s,Y_s) ds = \mu^*\lambda_1 = 0,
	\]
	which leads to
	\[
	\lim_{t \to \infty} \frac{1}{T} \int_0^T  X_s = K \left( 1 - \frac{\eps^2}{2} \right) - \int_0^{+ \infty} \frac{y}{1+x}d \mu^*(x,y).
	\]
	Now, due to the fact that $\mu^*(M_+)=1$, one has 
	\[
	\bar{y} :=  \int_0^{+ \infty} \frac{y}{1+x}d \mu^*(x,y) > 0
	\]
	and thus 
	\[
	\lim_{t \to \infty} \frac{1}{T} \int_0^T  (\hat X_s - X_s) ds = \bar y > 0.
	\]
	From this  we have
	\[
	\lim_{T \to \infty} \frac{1}{T} \int_0^T \left( \lambda_2(X_s, Y_s) - \lambda_2(\hat{X}_s, 0)\right) ds < 0.
	\]
	Indeed, let $C > 0$ such that 
	\[
	\lim_{T \to \infty} \frac{1}{T} \int_0^T  (\hat X_s - X_s) \1_{\{\hat X_s \leq C\}} ds  \geq \frac{\bar y}{2}.
	\]
	Then, it is easily seen that there exists $c > 0$ such that for all $0 \leq x \leq \hat x \leq C$, and all $y \geq 0$, one has $\lambda_2(x,y) - \lambda_2(\hat x,0) \leq -c( \hat x - x).$ In particular, by monocity of $H$ and the fact that $X_s \leq \hat X_s$ for all $s \geq 0$,  we have
	\begin{align*}
	\lim_{T \to \infty} \frac{1}{T} \int_0^T \left( \lambda_2(X_s, Y_s) - \lambda_2(\hat{X}_s, 0)\right) ds & \leq \lim_{T \to \infty} \frac{1}{T} \int_0^T \left( \lambda_2(X_s, Y_s) - \lambda_2(\hat{X}_s, 0)\right)\1_{\{\hat X_s \leq C\}} ds\\
	& \leq -\lim_{t \to \infty} \frac{c}{T} \int_0^T  (\hat X_s - X_s) \1_{\{\hat X_s \leq C\}} ds \\
	& \leq -\frac{c \bar y}{2}.
	\end{align*}
	We conclude that $\mu^* H=-\mu^*\lambda_2 \geq -\mu_x\lambda_2+\frac{c \bar y}{2}=\mu_x H + \frac{c \bar y}{2} > \mu_{\mathrm{x}} H$. This proves that when $\Lambda(\eps, K, \alpha) = 0$, $\Pcal_{inv}(\M_+)$ is empty. Moreover, we know that the process is persistent with respect to $\M_0^{\mathrm{x}}$. Putting this together, the only possible limit point for $(\Pi_t)_{t \geq 0}$ is $\mu_{\mathrm{x}}$. Furthermore, since the maps $(x,y) \mapsto \frac{(x+y)^2}{1+y}$ and $(x,y) \mapsto \frac{(x+y)^2}{1+x}$ are proper, Corollary \ref{cor:lambda0} and Remark \ref{rem:proper} imply that
	\[
	\lim_{T \to \infty} \frac{1}{T} \int_0^T Y_s ds =0
	\]
	and
	\[
	\lim_{T \to \infty} \frac{1}{T} \int_0^T X_s ds = \int_0^{+ \infty} x d \mu_{\mathrm{x}} (x) =  K \left( 1 - \frac{\eps^2}{2} \right).
	\]
	\QED
	\subsection{A stochastic model in a spatially heterogenous environments}
	\label{sub:patchy}
	In this section, we consider the example treated in \cite{HNY18} of a population submitted to random fluctuations of the environment and to spatio-temporal heterogeneity. The space is divided into $n$ patches, and the dynamics of the population within a patch follows a logistic SDE. There is also dispersal of the population, that is, individuals can move from one patch to the other. The precise model is the following. Let $X_t = (X_t^1, \ldots, X_t^n)$ be the vector of abundance in each patch at time $t$, then $X$ satisfy the SDE :
	\begin{equation}
	\label{eq:SDEpatch}
	d X_t^i =  \left[ X_t^i \left( a_i - b_i(X_t^i)\right) + \sum_{j=1}^n D_{j,i} X_t^j \right]dt + X_t^i dE_t^i,
	\end{equation}
	where $a_i > 0$ is the per capita growth rate in patch $i$, $b_i : \RR_+ \to \RR_+$ is the competition term in patch $i$, $D_{j,i} \geq 0$ is for $j \neq i$, the disperal rate of patch $j$ to patch $i$ and $E = \Gamma^T B$, where $\Gamma$ is a square $n \times n$ matrix and $B = (B^1, \ldots, B^n)$ is a standard Brownian motion. We also set $D_{i,i} = - \sum_{j \neq i} D_{j,i}$ and $\Sigma = \Gamma^T \Gamma$. 
	
	We work under the following assumptions, made in \cite{HNY18} :
	
	\begin{hypothesis}\
		\label{hyp:SDE}
		\begin{enumerate}
			\item For each $i \in \{1, \ldots, n\}$, $b_i : \RR_+ \to \RR_+$ is locally Lipschitz, vanishes only at $0$, and there exist constants $\gamma_b$ and $M_b$ such that, for all $x \in \RR_+^n$ with $\sum_i x_i \geq M_b$, one has
			\[
			\frac{\sum_{i=1}^n x_i (b_i(x_i) - a_i) }{\sum_{i=1}^n x_i} > \gamma_b;
			\]
			\item The matrix $D$ is irreducible;
			\item The matrix $\Sigma$ is non-singular.
		\end{enumerate}
	\end{hypothesis}
	
	These assumptions guarantee the existence of a unique strong solution to \eqref{eq:SDEpatch}, which moreover stays in $\RR_+^n$ if $X_0 \in \RR_+^n$.  As in \cite{HNY18}, we introduce the decomposition of the process : for any $x_0 \neq 0$ and $t \geq 0$, we set $S_t = \sum_i X_t^i$ and $Y_t^i = X_t^i/ S_t$. By Itô's formula, it can be shown that $(S_t, Y_t)$ evolves according to 
	\begin{equation}
	\label{eq:SDEpolar}
	\begin{cases}
	dY_t = \left[ \mathrm{Diag}(Y_t) - Y_t Y_t^T \right] \Gamma^T d B_t + D^T Y_t dt + \left[ \mathrm{Diag}(Y_t) - Y_t Y_t^T \right]\left( a - \Sigma Y_t - b(S_t Y_t)\right)dt\\
	dS_t = S_t \left( a  - b(S_t Y_t)\right)^T Y_t dt + S_t  Y_t^T \Gamma^T d B_t,
	\end{cases}
	\end{equation}
	where $Y_t = ( Y_t^1, \ldots, Y_t^n)$ lies in the simplex $$\Delta = \{ x = ( x_1, \ldots, x_n) \in \RR_+^n \: : x_1 + \ldots + x_n = 1 \},$$ and $ a:= (a_1, \ldots, a_n)$, $b(x):=(b_1(x), \ldots, b_n(x))$. It is now possible to extend equation \eqref{eq:SDEpolar} on $\{0\} \times \Delta$, by setting $S_t = 0$ and 
	\begin{equation}
	\label{eq:SDEboundary}
	d  Y_t = \left[ \mathrm{Diag}( Y_t) -  Y_t Y_t^T \right] \Gamma^T d B_t + D^T   Y_t dt + \left[ \mathrm{Diag}(  Y_t) -  Y_t Y_t^T \right]\left( a - \Sigma Y_t \right)dt.
	\end{equation}
	If we let $\tilde{X}_t$ be the solution to 
	\begin{equation}
	\label{eq:SDElinear}
	d \tilde X_t^i = \left[  a_i \tilde X_t^i + \sum_{j=1}^n D_{j,i} \tilde X_t^j \right]dt + \tilde X_t^i dE_t^i,
	\end{equation}
	and $\tilde{S}_t = \sum_i \tilde X_t^i$, then 
	\begin{equation}
	\label{eq:SDElinearS}
	d \tilde S_t = \tilde S_t  a^T \tilde Y_t dt + \tilde S_t  \tilde Y_t^T \Gamma^T d B_t,
	\end{equation}
	with $\tilde Y = Y$ subjected to \eqref{eq:SDEboundary}. It is proven in \cite{ERSS13} that $\tilde Y$ admits a unique invariant probabilty measure $\pi$ on $\Delta$. Set 
	\begin{equation}
	\label{eq:r}
	r = \int_{\Delta} \left( a^T y -  \frac{1}{2} y^T \Sigma^T y \right) d\pi(y). 
	\end{equation}
	In \cite{HNY18}, the authors show that the sign of $r$ determines the long term behaviour of $X$ : if $r < 0$, then the population abundance in each patch converges to $0$ exponentially fast, while if $r > 0$, the process $X$ admits a unique invariant probability measure on $\RR_{++}^n := \{ x \in \RR_+^n \: : x_i > 0\}$ and the law of $X$ converges polynomially fast to this stationnary distribution. The case $r = 0$ is not treated and left in the discussion as an open question.
	
	We show now that our method enables us to handle the critical case $r=0$:
	
	\bprop
	If $r= 0$, then, for all $i$, for all $x \in \RR_+^n$, $\PP_x$ - almost surely
	\[
	\lim_{t \to \infty} \frac{1}{t} \int_0^t X^i_s ds = 0.
	\]
	\eprop
	
	\prf First, let write the process in our background. We consider the process $(Z_t)_{t \geq 0} = ( S_t, Y_t)_{t \geq 0}$ defined on $\M = \RR_+ \times \Delta$, and evolving according to \eqref{eq:SDEpolar} on $\M_+ = \RR^*_+ \times \Delta$ and according to \eqref{eq:SDEboundary} on $\M_0 = \{0\} \times \Delta$. Proposition A.1 in \cite{HNY18} implies that Assumption \ref{hyp:feller} is satisfied under Assumptions \ref{hyp:SDE}. 
	
	One can check that for a function $f : \M \to \RR$, of class $C^2$ with bounded support, the generator $\L$ of $Z$ is given by :
	\[
	\L f(s,y) = \frac{\partial }{\partial s}f(s,y) s \left(a - b(sy)\right)^T y + \frac{1}{2} \frac{\partial^2 }{\partial s^2}f(s,y) s^2 y^T \Sigma y + A f(s,y),
	\]
	where $A f $ is a sum of terms, each of them involving at least one derivative of $f$ with respect to one of the coordinates of $y$. In particular, if $f(s,y) = g(s)$ for some function $g$, one has
	\[
	\L f(s,y) = g'(s) s \left(a - b(sy)\right)^T y +  \frac{1}{2} g''(s) s^2 y^T \Sigma y. 
	\]
	Let $\eps > 0$, and set $g(s) = (1+s)^{1+\eps}$ and $f(s,y) = g(s)$. Then, we get (formally) that
	\[
	\L f(s,y) = (1+\eps)f(s,y) \left[ \frac{s}{1+s}(a - b(sy))^T y + \left(\frac{s}{1+s}\right)^2\frac{1}{2} \eps y^T \Sigma y \right], 
	\]
	which by Assumption \ref{hyp:SDE} implies that
	\[
	\L f(s,y) \leq - \alpha f(s,y) + C,
	\]
	where $\alpha = \gamma_b - \frac{1}{2} \eps \| \Sigma \|$ is positive for $\eps$ small enough, and $C = \sup_{(s,y) \in [0,M_b] \times \Delta} \L f(s,y)$ is finite. From this, it is possible to prove that Assumption \ref{hyp:lyapinfty} is satisfied for $W(s,y) = (1+s)^{1+\eps}$, provided $\eps$ is small enough.

	Next, we prove that Assumption \ref{hyp:VH} is satisfied. We define two functions on $\M$:  
	\[
	H_1(s,y) = \frac{s}{1+s}  \left(a - b(sy)\right)^T y -  \frac{1}{2} \frac{s^2}{(1+s)^2}  y^T \Sigma y,
	\]
	and
	\[
	H_2(s,y) =  \left(a - b(sy)\right)^T y -  \frac{1}{2}  y^T \Sigma y.
	\]
	We define $V$ on $\M_+$ by setting $V(s,y) = \log (1+s) - \log s.$ By definition of $V$ and It\^o's formula,
	\[
	\L V(s,y) = H_1(s,y)-H_2(s,y).
	\]
	It is not hard to check that the functions $V$ and $H$ so defined satisfy Assumption \ref{hyp:VH}.
	We have from It\^o's formula that
	$$
	\lim_{T\to\infty}\Pi^z_T H_1=\lim_{T\to\infty}\frac{\EE_z \log (1+S_T)-\log(1+s)}T
	=0, z\in\M
	$$
	due to \cite[Lemma A.2]{HNY18}.
	As a result, $\nu H_1=0$ for any invariant probability measure $\nu$ on $\M$ of $(Z_t)_{t\geq0}$.
	Subsequently, we have, $r =-\pi H_2=  \pi H$, where $r$ is defined by \eqref{eq:r}, and by ergodicity of $\tilde{Y}$ and equation \eqref{eq:SDElinearS}, we have
	\begin{equation}
	r=\pi H = \lim_{ t \to \infty} \frac{1}{t} \int_0^t \left[ a^T \tilde Y_u - \tilde Y_u \Sigma \tilde Y_u \right] du = \lim_{ t \to \infty}  \frac{1}{t} \log \tilde{S}_t.
	\end{equation}
	Now we assume that $Z$ admits an ergodic invariant probability measure $\mu$ on $\M_+$. By the strong Feller property of $X$ on $\RR_{++}^n$, $\mu$ has to be unique, and thus the process is ergodic. In particular, we have
	\[
	-  \mu H=\mu H_2 = \lim_{ t \to \infty}  \frac{1}{t} \int_0^t \left[ \left( a - b( S_u Y_u)\right)^T Y_u - Y_u \Sigma Y_u \right] du. 
	\]
	Thus, to obtain the desired result that $\mu H > \pi H$,
	we will show that 
	\begin{equation}\label{eq:mugreaterpi}
	\begin{aligned}
	\lim_{ t \to \infty}& \frac{1}{t} \int_0^t \left[ a^T \tilde Y_u - \tilde Y_u \Sigma \tilde Y_u \right] du \\
	&>\lim_{ t \to \infty}  \frac{1}{t} \int_0^t \left[ \left( a - b( S_u Y_u)\right)^T Y_u - Y_u \Sigma Y_u \right] du.
	\end{aligned}
	\end{equation}
	While componentwise $a > a - b(S_u Y_u)$,
	\eqref{eq:mugreaterpi} is not straightforward because $\wdt Y_u - Y_u$ can be both negative and positive. The difficulty will be overcome by introducing an intermediate process to ease the comparison. For all $u \geq 0$, we set $\varsigma_u = \min_i b_i(S_u Y^i_u)$. Note that $\varsigma_u > 0$ by assumption on $b$. Now we introduce the process $\bar X = (\bar X^1, \ldots, \bar X^n)$ defined by
	\begin{equation}
	d \bar X_t^i =  \left[ \bar X_t^i \left( a_i - \varsigma_t  \right) + \sum_{j=1}^n D_{j,i} \bar X_t^j \right]dt + \bar X_t^i dE_t^i,
	\end{equation}
	
	By a classical comparison argument for SDE (see i.e. \cite{chueshov}) and positivity of $\varsigma_t$, we have $X_t^i \leq \bar X_t^i \leq \tilde{X}_t^i$ for all $t \geq 0$, provided the inequality holds at $0$. We also set $\bar{S}_t = \bar X_t^1 + \ldots + \bar X_t^n$, and then $S_t \leq \bar S_t \leq \tilde S_t$. Finally, we introduce $\bar Y = \bar X / \bar S$, which is well denifed as soon as $\bar X_0 \neq 0$. One can see that $\bar S$ and $\bar Y$ evolve according to
	\begin{equation}
	d \bar S_t =\bar  S_t \left( a  - \varsigma_t \right)^T \bar Y_t dt + \bar S_t  \bar Y_t^T \Gamma^T d B_t,
	\end{equation}
	\begin{equation}
	d \bar Y_t = \left[ \mathrm{Diag}(\bar Y_t) - \bar Y_t \bar Y_t^T \right] \Gamma^T d B_t + D^T \bar Y_t dt + \left[ \mathrm{Diag}(\bar Y_t) - \bar Y_t \bar Y_t^T \right]\left( a - \Sigma \bar Y_t - \varsigma_t \1 \right)dt,
	\end{equation}
	where $\1$ is the vector with all components equal to 1. Now, since $\bar Y_t \in \Delta$, one has $ (\mathrm{Diag}(\bar Y_t) - \bar Y_t Y_t^T )\1 =0$, thus
	\begin{equation}
	d \bar Y_t = \left[ \mathrm{Diag}(\bar Y_t) - \bar Y_t \bar Y_t^T \right] \Gamma^T d B_t + D^T \bar Y_t dt + \left[ \mathrm{Diag}(\bar Y_t) - \bar Y_t \bar Y_t^T \right]\left( a - \Sigma \bar Y_t  \right)dt,
	\end{equation}
	and by unicity of a strong solution to \eqref{eq:SDEboundary}, $\bar Y = \tilde Y$ almost surely whenever $\bar Y_0 = \tilde Y_0$. Thus we have
	\begin{align*}
	\lim_{t \to \infty} \frac{1}{t} \log(\bar S_t)& = \lim_{t \to \infty} \frac{1}{t} \int_0^t \left( a^T \tilde Y_u - \tilde Y_u \Sigma \tilde Y_u  - \varsigma_u\1^T \tilde{Y}_u \right) du\\ 
	& = - \pi H  - \lim_{t \to \infty} \frac{1}{t} \int_0^t \varsigma_u du\\
	& = - \pi H  - \int_{\M_+} \min_i b_i(s y_i) d \mu(y)\\
	& = - \pi H - \bar \varsigma,
	\end{align*}
	where $\bar \varsigma = \int_{\M_+} \min_i b_i(s y_i) d \mu(y) > 0$ because $\mu(\M_+) = 1$. Now, since $\bar S_t \geq S_t$, we have
	\begin{align*}
	- \pi H - \bar \varsigma & = \lim_{t \to \infty} \frac{1}{t} \log(\bar S_t)\\
	& \geq \lim_{t \to \infty} \frac{1}{t} \log( S_t)  = - \mu H,
	\end{align*}
	which yields $\mu H \geq \pi H + \bar \varsigma$.
	Thus, one can apply Corollary \ref{cor:lambda0} (and Remark \ref{rem:proper}): since the map $ (s,y) \mapsto \frac{(1+s)^{1+\eps}}{1+s}$ is proper, if $r=0$, one has  for all $(s,y) \in \M$, $\PP_{s,y}$ - almost surely
	\[
	\lim_{t \to \infty} \frac{1}{t }\int_0^t  S_u du = 0,
	\]
	or equivalently, for all $i$, for all $x \in \RR_+^n$, $\PP_x$ - almost surely
	\[
	\lim_{t \to \infty} \frac{1}{t} \int_0^t X^i_s ds = 0.
	\]
	\QED

	\subsection{SIS model in fluctuating environment}
	\label{sub:SIS}
	Here, we prove that the method used above also works in a SIS model with random switching environment. More precisely, we consider the model studied in \cite{BS17} and described as follows. Let $ d \geq 1$, $C = (C_{i,j})$ be an irreducible $d \times d$ matrix  with nonnegative entries and $D = (D_1, \ldots, D_d)$ a vector with positive entries. We define the vector field $F : \RR^d \to \RR^d$ by setting, for all $x \in \RR^d$, 
	\[
	F(x) = \left( C - \mathrm{Diag}(D) \right) x - \mathrm{Diag}(x) C x.
	\]
	This vector field was introduced by Lajmanovich and Yorke \cite{LajYorke} to describe a model of infection SIS (Susceptible - Infected - Susceptible) giving the evolution of a disease that does not confer immunity, in a population structured in $d$ groups. They analysed the differential equation on $[0,1]^d$ given by
	\[
	\dot x = F(x),
	\]
	that is,  componentwise,
	\[
	\frac{d x^i}{dt} = (1 - x^i) \left( \sum_{j=1}^d C_{i,j} x^j \right) - D_i x^i.
	\]
	In \cite{BS17}, we have considered a Piecewise Deterministic Markov Process $U = (X, \alpha)$ on $[0,1]^d \times \E$, where $\E = \{1, \ldots, N\}$ for some integer $N$ and evolving as follows :
	\begin{equation}
	\frac{d X_t}{dt} = F^{\alpha_t} ( X_t),
	\end{equation}
	where $\alpha$ is a Markov chain on $E$ and for all $k \in \E$, $F^k$ is the vector field defined like $F$ with $C$ and $D$ replaced by $C^k$ and $D^k$, respectively where $C^k$ and $D^k$ are a matrix and a vector as described above. We also set $A^k = C^k - \mathrm{Diag}(D^k)$. To analyse the long-term behaviour of $Z$, we have done in \cite{BS17} a polar decomposition : for $X_0 \neq 0$, we set $\rho_t = \| X_t \|$ and $ \Theta_t = \frac{X_t}{\rho_t}$. Then $W = (\rho, \Theta, \alpha)$ is still a PDMP, evolving according to 
	\beq
	\label{dThetarho}
	\left \{
	\begin{array}{l}
		\frac{d\Theta_t}{dt} = G^{\alpha_t}(\Theta_t) \\
		\frac{d \rho_t}{dt} =  \langle A^{\alpha_t} - \rho_t \mathrm{Diag}(\Theta_t) C   \Theta_t, \Theta_t \rangle \rho_t,
	\end{array}
	\right.
	\eeq
	where for
	all $i \in E$, $G^i$ is the vector field on $S^{d-1}$ defined by
	\beq
	\label{eq:defGi}
	G^i(\theta) = \left(A^i - \rho \mathrm{Diag}(\theta) C \right) \theta - \langle\left(A^i - \rho \mathrm{Diag}(\theta) C \right), \theta \rangle \theta.
	\eeq
	We set $\M_+ = \Psi( [0,1]^d \setminus \{0\} ) \times \E$, where $ \Psi : \RR^n \setminus \{0\} \to \RR_+^* \times S^{d-1}$ is defined by $\Psi(x) = ( \|x\|, \frac{x}{\|x\|})$. We also set $\M_0 = \{0 \} \times S^{d-1} \times \E$, then \eqref{dThetarho} can be defined on $\M_0$ be letting $\rho_t = 0$ for all $t \geq 0$ and 
	\begin{equation}
	\label{eq:polarlinear}
	\frac{d\Theta_t}{dt} =  A^{\alpha_t} \Theta_t -  \langle A^{\alpha_t}, \Theta_t, \Theta_t \rangle \Theta_t.
	\end{equation}
	We proved in \cite[Proposition 2.13]{BS17}, that on $\M_0 \simeq S^{d-1} \times \E$, the process $(\Theta, \alpha)$ admits a unique invariant probability $\pi$. We set 
	\[
	\Lambda = \int_{S^{d-1} \times \E} \langle A^i \theta, \theta \rangle d\pi(\theta, i).
	\]
	It has also be proven that the functions $V : \M_+ \to \RR_+$ and $H : \M \to \RR$, defined by $V(\rho, \theta, i)= - \log(\rho)$ and  by $H(\rho, \theta, i ) = - \langle A^i \theta, \theta \rangle + \rho \langle \mathrm{Diag}(\theta) C \theta, \theta \rangle$, respectively, satisfy asusmption \ref{hyp:VH}. It is easily seen that $\Lambda = - \pi H$. With our method, together with the results in \cite{BS17}, we can now fully describe the behaviour of $U$ according to the sign of $\Lambda$ :
	\bthm
	There are three possible asymptotic behaviours :
	\begin{enumerate}
		\item If $\Lambda < 0$, then for all $(x,i) \in [0,1]^d \times \E$, we have
		\[
		\PP_{x,i} \left( \limsup \frac{\log \|X_t\|}{t} \leq \Lambda\right) = 1.
		\]
		\item If $\Lambda = 0$, then for all $(x,i)\in [0,1]^d \times \E$, we have
		\[
		\lim_{t \to \infty} \frac{1}{t}\int_0^t \|X_s\| ds = 0 \quad \PP_{x,i} - \mbox{a.s.}
		\]
		and
		\[
		\PP_{x,i} - \lim_{t \to \infty} X_t = 0, 
		\]
		where $\PP_{x,i} - \lim$ denotes the convergence in probability. 
		\item If $\Lambda > 0$, then $U$ admits a unique invariant probability measure $\mu$ on $(0,1]^d \times E$. Moreover, there exists a Wasserstein distance $\mathcal{W}$ and $r > 0$ such that, for all probability $\nu$ with $\nu( \{0\} \times \E)$ and all $t \geq 0$, 
		\[
		\mathcal{W} ( \nu P_t, \mu) \leq e^{- rt} \mathcal{W}( \nu, \mu). 
		\]
	\end{enumerate}
	\ethm
	
	\prf
	The case $\Lambda < 0$ is Theoerem 4.3 in \cite{BS17}, while $\Lambda > 0$ is Theorem 4.12 in \cite{BS17}.
	
	To treat the case $\Lambda = 0$, we first prove that one can apply Proposition \ref{prop}. We assume that $W$ admits an invariant distribution $\mu$ on $\M_+$.   For all $t > 0$, we define 
	\[
	\varsigma_t = \min_{1 \leq i \leq d} X_t^i \left( \sum_j C_{i,j}^{\alpha_t} X_t ^j \right),
	\] 
	and we let $\bar X$ be the solution to
	\begin{equation}
	\frac{d \bar X_t }{dt } = \left( A^{\alpha_t} - \varsigma_t I \right) \bar X_t,
	\end{equation}
	where $I$ is the identity matrix of size $d$. We also let $Y$ be the solution to
	\begin{equation}
	\frac{d Y_t}{dt} = A^{\alpha_t} Y_t.
	\end{equation}
	By a comparison theorem for ordinary differential equations, we have $X_t^i \leq \bar{X}_t^i \leq Y_t$ for all $t \geq 0$, provided the inequality holds at time $0$. Finally, let $\bar \rho_t = \| \bar X_t\|$, $\bar \Theta_t = \frac{\bar X_t }{\bar \rho_t}$, $\tilde{\rho}_t = \| Y\|_t$ and $\tilde{\Theta_t} = \frac{Y_t}{\tilde \rho_t}.$ Then $\rho_t \leq \bar \rho_t \leq \tilde \rho_t$ and
	\[
	\frac{d \bar \Theta_t}{dt} = \left( A^{\alpha_t} - \varsigma_t \right) \bar \Theta_t - \langle \left( A^{\alpha_t} - \varsigma_t \right) \bar \Theta_t, \bar \Theta_t \rangle  \bar \Theta_t ,
	\]
	while $\tilde{\Theta}_t$ evolves according to  \eqref{eq:polarlinear}. Now, since $\langle \bar \Theta_t, \bar \Theta_t \rangle = 1$ for all $t \geq 0$, we can see that $\bar \Theta_t$ is also driven by \eqref{eq:polarlinear}, thus $\bar \Theta_t = \tilde \Theta_t$ for all $t \geq 0$ whenever $\bar \Theta_0 = \tilde \Theta_0$. On the other hand, one can check that 
	\[
	\lim_{t \to \infty} \frac{\log \bar \rho_t}{t} = \lim_{t \to \infty} \frac{1}{t} \int_0^t \langle \left( A^{\alpha_s} - \varsigma_t \right) \bar \Theta_s, \bar \Theta_s \rangle ds.
	\]
	We also have
	\[
	- \pi H = \lim_{t \to \infty} \frac{\log \tilde \rho_t}{t} = \lim_{t \to \infty} \frac{1}{t} \int_0^t \langle A^{\alpha_s} \tilde \Theta_s, \tilde \Theta_s \rangle ds,
	\]
	Without loss of generality, one may assume that $\mu$ is ergodic, and therefore, one has for $\mu$ almost every $(\rho_0, \theta_0, i) \in \M_+$, $\PP_{(\rho_0, \theta_0, i)}$ - almost surely,
	\[
	\lim_{t \to \infty} \frac{1}{t} \int_0^t \varsigma_t = \int_{\M_+} \rho^2 \min_i \theta^i \left( \sum_j C_{i,j}^k \theta^j\right) d \mu(\rho, \theta, k) : = \bar \varsigma.
	\]
	Then, $\bar \varsigma > 0$ because on $\M_+$, $\rho > 0$ and $\mu(\{(\rho, \theta, i) \in \M_+ \: : \theta^i > 0\}) = 1$ since $\partial S^{d-1}$ is transient for $W$. Thus, due to the fact that $\bar \Theta_t = \tilde \Theta_t$, we get for $\mu$ almost every   $(\rho_0, \theta_0, i) \in \M_+$, $\PP_{(\rho_0, \theta_0, i)}$ - almost surely,
	\[
	\lim_{t \to \infty} \frac{\log \bar \rho_t}{t} = - \pi H  - \bar \varsigma,
	\]
	which combined with 
	\[
	\lim_{t \to \infty} \frac{\log  \rho_t}{t} = - \mu H \quad \PP_{(\rho_0, \theta_0, i)} - \mbox{a.s.}
	\]
	and $\bar \rho_t \geq \rho_t$ gives $ \mu H \geq \pi H + \varsigma > \mu H$. Thus, by Proposition \ref{prop}, $\Lambda > 0$. Hence, if $\Lambda = 0$, the unique stationnary distribution of $W$ is $\pi$, which is concentrated on $\M_0$. In particular, going back to the process $U$, its unique invariant distribution is $\delta_0 \otimes p$, where $p$ is the unique stationnary distribution of $\alpha$ on $\E$. In particular, for all $(x,i) \in [0,1]^d \times \E$, one has $\PP_{x,i}$ - almost surely that 
	\begin{equation}
	\label{eq:cvcesaro}
	\lim_{t \to \infty} \frac{1}{t} \int_0^t \| X_s \| ds = 0.
	\end{equation}
	
	To prove that $X$ converges in probability to $0$, we use results on monotone random dynamical systems due to Chueshov \cite{chueshov}. Let $\Omega = \mathbb{D}( \RR_+, \E)$ be the Skorhokhod space of càdlàg functions $ \omega : \RR_+ \to \E$, endowed with its Borel sigma field $\Fcal$, and on which we define the shift $\mathbf{\Theta} = (  \mathbf{\Theta}_t)_{t \geq 0}$ by  
	\[
	\mathbf{\Theta}_t ( \omega)(s) = \omega(t+s).
	\]
	We let $\PP_p$ be a probability measure on $(\Omega, \Fcal)$ such that the canonical process $I$ has the law of $\alpha$ starting from its ergodic probability measure $p$. Then, the process $\Psi(\omega, t)$ defined by
	\begin{equation}
	\begin{cases}
	\frac{d \Psi(t,\omega) x}{dt} = F^{\omega(t)}( \Psi(t,\omega) x)\\
	\Psi(0,\omega) x = x
	\end{cases}
	\end{equation}
	is a Random Dynamical System over the ergodic dynamical system $(\Omega, \Fcal, \PP_p, \mathbf{\Theta})$ (see e.g. \cite{ards} for definitions and the thesis of the secound author \cite[Section 1.4]{these} for more details on random dynamical systems and links with PDMPs). Moreover, the proporties of $F$ make $\Psi$ a monotone subhogeneous random dynamical system (see \cite[Section 4]{BS17}) for which $\1 =(1, \ldots, 1)$ is a super-equilibrium. That is, $\Psi(t, \omega) \1 \leq \1$ for all $t \geq 0$ and $\omega \in \Omega$ (see \cite[Definition 3.4.1]{chueshov}). Moreover, for all $t \geq 0$ and all $\omega \in \Omega$, 
	\[
	\Psi(t, \omega) \left( [0,1]^d \setminus \{0\} \right) \subset (0,1)^d.
	\]
	Hence, it is easily to check that we can apply Proposition 5.5.1 in \cite{chueshov}. According to this result, either, for all $x \in [0,1]^d$, 
	\beq
	\label{eq:cv0}
	\lim_{ t \to \infty} \Psi(t, \mathbf{\Theta}_{-t} \omega) x = 0
	\eeq
	or, there exists an equilibrium $u(\omega) \gg 0$ such that, for all $x > 0$ and all $\omega \in \Omega$, 
	\begin{equation}
	\label{eq:cvu}
	\lim_{ t \to \infty} \Psi(t, \mathbf{\Theta}_{-t} \omega) x = u( \omega).
	\end{equation}
	Now, assume that \eqref{eq:cvu} holds. In particular, by dominated convergence and invariance of $\PP_p$ under $\mathbf{\Theta}$, one has on the one hand
	\begin{equation}
	\label{eq:cvEpsi}
	\lim_{t \to \infty} \EE_p\left( \| \Psi(t, \omega) x \|  \right) = \EE_p(\|u\|) > 0.  
	\end{equation}
	On the other hand, one can check that the law of $X_t$ under $\PP_{x,p}$ is the same as the law of $\Psi(t, \cdot)x$ under $\PP_p$. In particular, 
	\begin{equation}
	\label{eq:egalite}
	\EE_{(x,p)}\left( \| X_t \| \right) = \EE_p \left( \| \Psi(t, \omega) x \| \right). 
	\end{equation}
	Thus, \eqref{eq:cvEpsi} and \eqref{eq:egalite} imply that
	\[
	\lim_{t \to \infty} \frac{1}{t} \int_0^t \EE_{(x,p)}\left( \| X_s \| \right) ds = \EE_p(\|u\|)>0,
	\]
	which is in contradiction (by dominated convergence)  with \eqref{eq:cvcesaro}. Hence, \eqref{eq:cv0} holds. This and \eqref{eq:egalite} yield that for all continuous map $f : [0,1]^d \to \RR$,
	\[
	\lim_{t \to \infty} \EE_{x,p}\left( f(X_t) \right)  = f(0), 
	\]
	which implies that $X_t$ converges in law, hence in probability, to $0$, under $\PP_{x,p}$. It is easily seen that one can now replace $p$ by any starting point $i$. \QED
	\subsection{SEIR model with switching}
	\label{sub:SEIR}
	SEIR models describe the dynamics of an infectious disease
	with which individuals experience a long incubation duration (the “exposed” compartment). The classical SEIR model
	consists of the following differential equations for 4 classes of individuals (Susceptible - Exposed - Infectious - Recovered):
	
	\begin{equation}
	\begin{cases}
	\dot S=&\Lambda- \gamma S  - \beta S I\\
	\dot E=& \beta S I-(\gamma+\delta)E\\
	\dot I=& \delta E -(\gamma+\gamma_1) I\\
	\dot R=& \gamma_1 I- \gamma R
	\end{cases}
	\end{equation} 
	where $\Lambda,\gamma, \beta, \delta,\gamma_1$ are positive constant. We refer to \cite{schwartz1983infinite,hethcote2000mathematics}
	for details about this model and its variants. 
	In contrast to stochastic SIR and SIRS models, which have been studied extensively,
	few papers deal with 
	stochastic SEIR models because
	standard arguments used to treat SIR and SIRS models do not seem effective for SEIR models.
	In this section,
	we wish to consider an SEIR model in a switching environment.
	Let $N$ be a positive integer, and set $\E=\{1,\ldots,N\}$. Let
	$(\alpha_t)_{t \geq 0}$ be a irreducible Markov chain on $\E$ and consider the following system
	
	\begin{equation}\label{e:SEIR}
	\begin{cases}
	\dot S=&\Lambda- \gamma S  - \beta(\alpha_t) S I\\
	\dot E=& \beta(\alpha_t) S I-(\gamma+\delta(\alpha_t))E\\
	\dot I=& \delta(\alpha_t) E -(\gamma +\gamma_1(\alpha_t))I
	\end{cases}
	\end{equation} 
	where the component $R$ is removed because it does not affect the dynamics of the others.
	
	Let $U_t=E_t+I_t$ and $V_t=\frac{I_t}{U_t}$, $Z_t=(S_t, V_t, U_t,\alpha_t)$, we can rewrite \eqref{e:SEIR} as
	
	\begin{equation}\label{e:SEIR}
	\begin{cases}
	\dot S=&f_S(Z_t)\\
	\dot V=& f_V(Z_t)\\
	\dot U=& U_tf_U(Z_t)
	\end{cases}
	\end{equation} 
	where 
	$$f_S(z)=\Lambda- \gamma s - \beta(k) su(1-v),$$
	$$f_U(z)=(\beta(k)s-\gamma_1(k)-\gamma)v-\gamma (1-v)=(\beta(k)s-\gamma_1(k))v -\gamma $$
	$$f_V(z)=(\sigma(k)(1-v)-\gamma v-\gamma_1(k) v)- vf_U(z)=\sigma(k)(1-v)-\gamma_1(k)v -(\beta(k)s-\gamma_1(k))v^2,$$
	and $z=(s, u, v, k)$.
	For this system, we have
	\begin{equation}
	\label{eM.5}
	\M:=\left\{z\in\RR^2_+\times[0,1]\times\E: s+u\leq \frac\Lambda\gamma\right\}\, 
	\text{ and } \M_0:=\{z\in\M: u=0\}.
	\end{equation}
	
	In this model, 
	$H(z):=-f_U(z)$ and
	$\mathcal V(z):=\log\frac\Lambda\gamma-\log u$
	satisfy Assumption \ref{hyp:VH}.
	Unlike the arguments in Subsections 3.3 and 3.4,
	it does not seem practically possible to treat the critical case by introducing an intermediate process.
	Because the function $f_U(z)$ is increasing in $s$ while
	$f_V(z)$ is decreasing in $s$, we introduce the following function:
	\beq\label{wH}
	\wdt H(z)=-f_U(z)-\frac{f_V(z)}v= -\sigma(k)\frac{1-v}v+\gamma+\gamma_1(k).
	\eeq
	If $U_0=0$ then $U_t=0, t\geq0$ and $\lim_{t\to\infty} S_t=\frac\Lambda\gamma$.
	Let $\wdt V_t$ be the solution to
	$$
	\dot{\wdt{V}}=f_V\left(\frac{\Lambda}{\gamma}, 0,\wdt V_t,  \alpha_t\right).$$
	Then, one can show that $(\wdt V_t, \alpha_t)$ has a unique invariant measure $\pi_V$ on $[0,1]\times\E$
	(see e.g. \cite[Proposition 2.1]{BL16} or \cite{DN11}).
	Moreover, since $f_V(z)=\sigma(k)>0$ if $z=(s,u,v,k)$ with $v=0$,
	there exists $v_0>0$ such that
	$\liminf_{t\to\infty} V_t\geq v_0>0$
	for any initial value $z\in\M.$
	As a result,
	$$\PP_z\left\{\lim_{T\to\infty}\frac1T\int_0^T \frac{f_V(Z_t)}{V_t}=\lim_{T\to\infty}\frac{\log V_t}{T}=0\right\}, z\in\M.$$
	Hence, for any invariant probability measure $\mu$ of $(Z_t)_{t \geq 0}$, we have
	\beq\label{wH2}\int_\M \frac{f_V(z)}v \mu(dz)=0, \text{ or equivalently } \mu H=\mu\wdt H.
	\eeq
	
	Then
	$\pi:=\bdelta_{(\frac\Lambda\gamma,0)}\otimes\pi_V$ is the unique invariant measure on $\M_0$.
	By the ergodicity of $(\wdt V_t, \alpha_t)$ and \eqref{wH},\eqref{wH2}
	we have
	\beq\label{4-e3}
	\begin{aligned}
		\wdt\Lambda:=-\pi H=&-\pi\lim_{T\to\infty}\int_0^T f_U\left(\frac{\Lambda}{\gamma},0, \wdt V_t,  \alpha_t\right)dt\\
		=&
		\lim_{T\to\infty}\int_0^T\left(\sigma(\alpha_t)\frac{1- \wdt V_t}{\wdt V_t}-\gamma -\gamma_1(\alpha_t)\right)dt.
	\end{aligned}
	\eeq
	
	With $\M, \M_0$ define in \eqref{eM.5}, we have the following theorem
	\bthm
	\begin{enumerate}
		\item If $\wdt\Lambda < 0$, then for all $0 < \lambda < - \wdt \Lambda$, there exist $\eta > 0$ and $r > 0$ such that, for all  $z\in\M_+:=\M\setminus\M_0$ with $u \leq r$,  we have
		\[
		\PP_{x,i} \left( \limsup \frac{\log U_t}{t} \leq  - \lambda\right) \geq \eta.
		\]
		\item If $\wdt\Lambda = 0$, then for all $z\in\M_+\times\E$, we have
		\[
		\lim_{T \to \infty} \frac{1}{T}\int_0^T U_tdt = 0 \quad \PP_{z} - \mbox{a.s.}
		\]
		\item If $\wdt\Lambda > 0$, then $Z$ is $H$-persistent and it admits an invariant probability measure on $\M_+$.
	\end{enumerate}
	\ethm
	\prf
	We start by proving the first and third claims. For $\alpha \in \E$, we define the vector field
	\[
	F^{\alpha}(s,e,i) =
	\begin{cases}
	\Lambda - \gamma s + \beta(\alpha) s i\\
	\beta(\alpha) s i - ( \gamma + \delta(\alpha) ) e\\
	\delta(\alpha) e - ( \gamma + \gamma_1(\alpha)) i
	\end{cases} 
	\]
	Then, letting $(X_t)_{t \geq 0} = (S_t, E_t, I_t)_{t \geq 0}$, we have $ \dot{X}_t = F^{\alpha_t}(X_t)$. 
	Note that $(\frac{\Lambda}{\gamma},0,0)$ is a common equilibrium of the vector fields $F^{\alpha}$ and that the line $\RR_+ \times \{(0,0)\}$ is invariant for each of the vector fields. This is exactly the setting of application of the results in \cite{S18}. The Jacobian matrix of $F^{\alpha}$ at $(\frac{\Lambda}{\gamma},0,0)$ is given by
	\[
	A^{\alpha} = \begin{pmatrix}
	- \gamma & - \beta(\alpha)\frac{\Lambda}{\gamma} & 0\\
	0 & -( \gamma + \delta(\alpha) ) & \beta(\alpha)\frac{\Lambda}{\gamma}\\
	0 & \delta(\alpha) & ( \gamma + \gamma_1(\alpha))
	\end{pmatrix}.
	\]
	We let $D = (\gamma)$, $C^{\alpha} = (- \beta(\alpha)\frac{\Lambda}{\gamma}, 0)$ and 
	\[
	B^{\alpha} = \begin{pmatrix}
	-( \gamma + \delta(\alpha) ) & \beta(\alpha)\frac{\Lambda}{\gamma}\\
	\delta(\alpha) & ( \gamma + \gamma_1(\alpha))
	\end{pmatrix},
	\]
	so that 
	\[
	A^{\alpha} = \begin{pmatrix}
	D & C^{\alpha}\\
	0 & B^{\alpha}
	\end{pmatrix}.
	\]
	Finally, we define $ \Lambda_D = - \gamma$ and $\Lambda_B = \int \langle B^{\alpha} \theta, \theta \rangle d\pi( \alpha, \theta)$, where $\pi$ is the unique invariant probability measure of the process $(\Theta, \alpha)$, where $\Theta$ is subjected to \eqref{eq:polarlinear} with $A^{\alpha}$ replaced by $B^{\alpha}$ ( the uniqueness of $\pi$ comes from the particular form of $B$, see \cite[Proposition 2.13]{BS17}). Then, $\tilde \Lambda = \Lambda_B$. Indeed, $\tilde \Lambda$ is defined as the growth rate of $U$, which is the $L^1$-norm of $(E,I)$, while $\Lambda_B$ is defined as the  growth rate of $U_2$, the $L^2$-norm of $(E,I)$. By equivalence of the norm on $\RR^2$, we must have $ \tilde{\Lambda} = \Lambda_B$. The third claim is hence a direct application of Theorem 2.8 in \cite{S18}. The first claim follows from Theorem 2.7 in \cite{S18}.

	Now, we prove the second claim. Let  assume that $\Pcal_{inv}(\M_+)$ is nonempty with an ergodic measure $\mu$.
	
	Since
	$f_V(z)+\beta(k)\left(\frac\Lambda\gamma-s)\right)v^2= f_V(\frac\Lambda\gamma, 0, v, k)$ for any $z=(s, u, v, k)\in\M$,
	we have
	$V_t\geq \wdt V_t$ given $V_0\geq \wdt V_0$.
	Let $Z_t$ have the intial distribution $\mu$
	and $\wdt V_0=V_0$.
	
	By the ergodicity we have
	$$
	\lim_{T\to\infty}\frac1T\int_0^T(V_t-\wdt V_t)dt=\int_\M v\mu(dz)-\int_\M v\pi(dz) \text{ a.s.}.
	$$
	We will show that $\int_\M v\mu(dz)-\int_\M v\pi(dz)>0$
	by a contradiction argument.
	Note that, since
	$$
	|f_V(s, u, v, k)-f_V(s, u, \wdt v, k)|\leq C|v-\wdt v|
	$$
	for some constant $C>0$,
	we have
	\beq\label{e4-13}
	\limsup_{T\to\infty}\frac1T\int_0^T\left|f_V(Z_t)dt - f_V\left(S_t, U_t, \wdt V_t, \alpha_t\right)\right|dt
	\leq 
	\lim_{T\to\infty}\frac1T\int_0^T C\left(V_t-\wdt V_t\right)=0
	\eeq
	if $\int_\M v\mu(dz)-\int_\M v\pi(dz)=0$.
	On the other hand,
	$$f_V(s, u, \wdt v, k)=f_V\left(\frac{\Lambda}\gamma, 0, \wdt v, k\right)
	+\left(\frac{\Lambda}\gamma-s\right)\beta(k)\wdt v^2
	$$
	which leads to
	\beq\label{4-e12}
	\begin{aligned}
		\lim_{T\to\infty}\frac1T\int_0^T f_V\left(S_t, U_t, \wdt V_t, \alpha_t\right) dt
		=&
		\lim_{T\to\infty}\frac1T\int_0^T f_V\left(\frac\Lambda\gamma, 0, \wdt V_t, \alpha_t\right) dt\\
		&+\lim_{T\to\infty}\frac1T\int_0^T\left(\frac{\Lambda}\gamma-S_t\right)\beta(\alpha_t)\wdt V_t\\
		>&\lim_{T\to\infty}\frac1T\int_0^T f_V\left(\frac\Lambda\gamma, 0, \wdt V^2_t, \alpha_t\right) dt.
	\end{aligned}
	\eeq
	where we  use the ergodicity of $(Z_t,\wdt V_t)$ on $\M_+\times(0,1)$
	to have that
	$$
	\lim_{T\to\infty}\frac1T\int_0^T\left(\frac{\Lambda}\gamma-S_t\right)\beta(\alpha_t)\wdt V_t>0$$
	Combining \eqref{e4-13} and \eqref{4-e12} we have
	$$
	\lim_{T\to\infty}\frac1T\int_0^Tf_V(Z_t)dt - \lim_{T\to\infty}\frac1T\int_0^Tf_V\left(\frac\Lambda\gamma, 0, \wdt V_t, \alpha_t\right)dt>0
	$$
	if $\int_\M v\mu(dz)-\int_\M v\pi(dz)=0$.
	
	However, it contradicts the fact that
	$$
	\begin{aligned}
	\lim_{T\to\infty}\frac1T\int_0^Tf_V(Z_t)dt -& \lim_{T\to\infty}\frac1T\int_0^Tf_V\left(\frac\Lambda\gamma, 0, \wdt V_t, \alpha_t\right)dt\\
	&=\int_\M f_V(z)\mu(dz)-\int_\M f_V(z)\pi(dz)=0-0
	\end{aligned}
	$$
	where the last equality is due to an argument similar to \eqref{wH2}.
	Thus,
	\beq\label{4-e10}
	\lim_{T\to\infty}\frac1T\int_0^T(V_t-\wdt V_t)dt>0 \text{ a.s.}.
	\eeq
	Since $\wdt H$ is an increasing function in $v$ with positive derivative,
	we can easily implies from \eqref{4-e10} and the fact that $V_t\geq \wdt V_t$ that
	$$
	\mu\wdt H-\pi\wdt H=\lim_{T\to\infty}\frac1T\int_0^T\wdt H(Z_t)dt - \lim_{T\to\infty}\frac1T\int_0^T\wdt H\left(\frac\Lambda\gamma, 0, \wdt V_t, \alpha_t\right)dt>0
	$$
	%
	%
	%
	%
	%
	%
	%
	In view of Corollary \ref{cor:lambda0}, 
	we obtain the second claim of the theorem.
	The proof is complete.
	\QED
	
	\section{Conclusion}	
	
	In this paper, we have given a general method to deal with the critical case in population dynamics in random environment. We apply the method to five different models, including epidemiological, prey-predator, and population in structured environment. 
	
	When our results apply, there is extinction in temporal average in the critical case. A natural question is wether it is possible to find other results, such that there is persistence (maybe in a weaker sense) in the critical case. 
	
	Our method consists in looking at integrals of the function $H = \L V$ with respect to invariant measures of the process. For some models, another method is possible, as used for example for some PDMP in \cite{HK19}. The idea is the following. Assume that if $\Pcal_{inv}(\M_+)$ is nonempty, then it is possible to compute, or at least, to estimate, the density of an invariant probability $\mu \in \Pcal_{inv}(\M_+)$. Then, this density must satisfy some integrability conditions, which can be violated if $\Lambda^+(H) = 0$ (see e.g \cite[Theorem 3.1]{HK19} or \cite[Lemma 6]{GPS19}). Hence, if  $\Lambda^+(H) = 0$, $\Pcal_{inv}(\M_+)$ has to be empty. This alternative method is close in spirit to ours, since it comes to a condraction when assuming that $\Pcal_{inv}(\M_+)$ is nonempty and $\Lambda^+(H) = 0$.

	\bibliographystyle{amsalpha}
	\bibliography{bibliothese}
\end{document}